\documentclass[12pt]{amsart}
\usepackage{amssymb}
\usepackage{amsfonts}
\usepackage{amsmath}
\usepackage{graphicx}

\newtheorem{theorem}{Theorem}[section]

\newtheorem{step}{}[section]
\newtheorem{claim}[theorem]{Claim}

\newtheorem{conjecture}[theorem]{Conjecture}
\newtheorem{corollary}[theorem]{Corollary}

\newtheorem{lemma}[theorem]{Lemma}

\newtheorem*{problem*}{Problem}

\newcounter{ctr}

\def\be{\begin{equation}}
\def\ee{\end{equation}}
\def\QATOP#1#2{{#1 \atop #2}}

\begin{document}
\author{Jonah Blasiak}
\title{A special case of Hadwiger's conjecture}
\maketitle

\begin{abstract}
We investigate Hadwiger's conjecture for graphs with no stable set of size
3. Such a graph on at least $2t-1$ vertices is not $t-1$ colorable, so is
conjectured to have a $K_t$ minor. There is a strengthening of Hadwiger's conjecture in this case, which states that there is always a minor in which the
preimages of the vertices of $K_t$ are connected subgraphs of size one or
two. We prove this strengthened version for graphs whose complement has an
even number of vertices and fractional chromatic number less than 3. We
investigate several possible generalizations and obtain counterexamples for
some and improved results from others. We also show that for sufficiently large $n=|V(G)|$,
a graph with no stable set of size 3 has a $K_{\frac{1}{9} \; n^{4/5}}$ minor using only sets of size one or two as preimages of vertices.
\end{abstract}
\maketitle

\section{Introduction}
A graph $H$ is a minor of a graph $G$ if $H$ can be obtained from a subgraph of $G$ by contracting edges.  If $E$ is the set of contracted edges, we call a connected component of the graph $(V(G),E)$ a \emph{prevertex}. Upon contraction of $E$, each prevertex becomes a vertex of $H$. It is clear from the many excluded minor theorems, the connections between minors and surface embeddings, and Robertson and Seymour's Graph Minor Theorem (see e.g.\cite{M},\cite{T}) that studying minors is an excellent way to study graph structure.  Perhaps the first important result to make use of minors was Kuratowski's theorem.  We state Wagner's reformulation (see e.g. \cite{T})
\begin{theorem}
A graph is planar if and only if it has no $K_5$ or $K_{3,3}$ minor.
\end{theorem}

Long before this, in 1852, Francis Guthrie formulated the four-color theorem, which states that every loopless planar graph is 4-colorable \cite{T}. Kuratowski's theorem gave a new way to study the four-color conjecture. In 1937 Wagner proved that the statement ``every loopless graph with no $K_5$ minor is 4-colorable'' is equivalent to the four-color conjecture (see e.g. \cite{RST}).  In 1943 Hadwiger and Dirac proved that every loopless graph with no $K_4$ minor is 3-colorable, and Hadwiger conjectured that (see e.g. \cite{RST})

\begin{conjecture}
For $t \geq 1$, every loopless graph with no $K_{t}$ minor is $(t-1)$-colorable.
\end{conjecture}

If true, this is a marvellously simple connection between complete minors and chromatic number, and is therefore considered one of the most important problems in graph theory. It has proved to be as difficult as it is beautiful.
In 1977, Appel and Haken proved the four-color conjecture, but it was extremely complicated and computer assisted (see e.g. \cite{RST}, \cite{T}). Robertson, Sanders, Seymour, and Thomas cleaned up the proof, but it is still computer assisted \cite{RSST}.  In 1993, Robertson, Seymour, and Thomas proved that the four-color theorem is equivalent to Hadwiger's conjecture for $t=5$ \cite{RST}. At present, Hadwiger's conjecture has been proved for $t\leq 5$ and is open for all $t \geq 6$.

We investigate Hadwiger's conjecture in another regime---when $t$ is comparable to the number of vertices in the graph. It is thought that if Hadwiger's conjecture is false, this is the most likely place for a counterexample. We restrict our attention to the case where $G=(V,E)$ has no stable set of size 3.  This implies that there are at least $|V|/2$ color classes in a proper coloring of $G$; Hadwiger's conjecture implies that $G$ has a complete minor of size at least $|V|/2$. A strengthening conjectured by Seymour is

\begin{conjecture}
\label{mainconjecture} If $G=(V,E)$ has no stable set of size 3, then $G$ has a complete minor of size at least $|V|/2$ using only edges or single vertices as prevertices.
\end{conjecture}

We call this the SSH conjecture; SS stands for Seymour's strengthening and stable set and H stands for Hadwiger. Our main result states that SSH is true if the edges of $G$ can be partitioned into two sets with certain properties. We also show that this condition is satisfied by some reasonably interesting classes of graphs (graphs whose complement is 3-colorable, for example).

We strongly believe SSH is true because many attempts at constructing counterexamples have failed.
However, our intuition for graphs with no stable set of size three is
severely limited. We have much difficulty constructing graphs that our results do not apply to---graphs with no stable set of size 3, large connectivity, no dominating edges, and large chromatic number in the complement. The only random graphs we can construct with these properties are extremely dense and have large complete minors.

Before stating our results, we need some notation.  If $A$ and $B$ are sets, $A$
\emph{intersects} $B$ means $A \cap B \neq \emptyset$. $[n]$ will denote the set $
\{1,2,\ldots n\}$.

All graphs in this thesis are finite. Let $G$ be a graph. We will sometimes write $G = (V,E)
$, which means $G$ has vertex set $V$ and edge set $E$; we will
also use $V(G)$ and $E(G)$ for the vertex and edge sets of $G$.
When there is no ambiguity, we use $n$ instead of $|V|$ without saying so explicitly.
If $S \subseteq V(G)$, $G[S]$ is the induced subgraph $G \backslash (V(G)-S)$. $\overline{G}$ is the complement of $G$. $d_G(v)$ is the degree of $v$ in $G$, and the subscript $G$ will be omitted when there is no ambiguity.  We will write
$(u,v)$ for an edge with ends $u$ and $v$, and $u \sim v$ ($u \nsim v$) means edge $(u,v)$ is (is not) present. If $U$ and $V$ are disjoint
vertex sets, a $(U,V)$ edge is some edge with one end in $U$ and one end in
$V$; the $(U,V)$ edges is the set of all edges with one end in $U$ and one
end in $V$.

We will say that the vertex sets $U$ and
$V$ \emph{touch} if they intersect or there is some edge with an end in each set.
We will also speak of two edges touching or an edge and a vertex touching;
we just identify the edge $(u,v)$ with the set $\{u,v\}$ and use the notion
of touching just mentioned. We say $U$ is \emph{complete} (\emph{anticomplete}) to $V$
if every (no) edge $(u,v)$ $u \in U, v \in V$ is present. If $v$
is a vertex, $N(v)$ will denote its set of neighbors (and will not include $v
$); if $V$ is a vertex set, $N(V) = \bigcup_{v \in V} N(v)$. A \emph{dominating edge} of $G$ is an edge that touches every vertex of $G$.

Vertices $u$, $v$ are said to be \emph{twins} if they are non-adjacent and $N(u)=N(v)$. $G^{\prime}$ is
a \emph{blown up} $G$ if $G$ can be obtained from $G^{\prime}$ by
identifying pairs of twin vertices. Vertex \emph{duplication} is the action of
replacing a vertex by two non-adjacent vertices with the same neighbors as
the original. Unfortunately, these are the standard definitions of twins and duplication, but we want a ``complementary'' definition. We say vertices $u$, $v$ are \emph{c-twins} if they are adjacent and $N(u)=N(v)$; c stands for complement and clique.  Define \emph{c-duplication} and \emph{c-blown up} similarly.

An \emph{antitriangle} is a stable set of size 3. Let $\mathfrak{A}$ be the set of graphs with no antitriangle.

\section{First observations}
A simple but important observation is that a minimal counterexample to SSH has no dominating edges.
In fact, we can win in two ways. If $G$
has a dominating edge, $e= (u,v)$, then we can use $e$ as a prevertex
together with a minor on $G \backslash \{u, v\}$ found inductively. Or we
observe that $G \backslash e \in \mathfrak{A}$, and by induction find
a complete minor on it.

Another preliminary result gives a lower bound on the connectivity of a
counterexample to SSH.

\begin{lemma}
\label{cutset} If $G=(V,E) \in \mathfrak{A}$ and has a cut set, $M$, of size at most $\frac{n}{2}$, then SSH holds.
\end{lemma}

\begin{proof}
Choose $M$ as small as possible. Let $L,R$ be a partition of $V-M$ such that $%
L $ and $R$ don't touch. $G$ has no antitriangle implies that $L$ and $R$
are cliques and that every vertex in $M$ is either complete to $L$ or complete
to $R $. Let $M_L,M_R$ partition $M$ so that every vertex in $M_{L}$ ($M_{R}$) is complete to $L$ ($R$).
Any $A\subseteq M_{L}$ of size at most $|R|$ is matchable into $R$. If not,
by Hall's matching condition, $\exists S\subseteq A$ such that $%
|S|>|N(S)\cap R| $. But then $(M-S)\cup |N(S)\cap R|$ is a cutset
because it separates $L\cup S$ and $R-N(S)$; it is smaller than M, contradiction. Now let $Y$ be a matching from $%
M_{L}$ to $R$ of size $\min ({|M_{L}|,|R|)}$. The vertices of $L$, together
with the edges of $Y$ are the prevertices of a complete minor (any pair of
edges in $Y $ is adjacent because they both have an end in the clique $R$).
We can, of course, do the same thing with vertices from $R$ and a matching
from $M_{R}$ to $L$. So without loss of generality $|L|+|M_{L}|\geq
|R|+|M_{R}|$. The size of the complete minor is $|L|+\min ({|M_{L}|,|R|)}%
=\min ({|L|+|M_{L}|,|L|+|R|)}\geq \frac{n}{2}$ by the assumption that $%
|M|\leq \frac{n}{2}$.
\end{proof}

\begin{figure}[tbp]
\begin{center}
\includegraphics [scale=1]{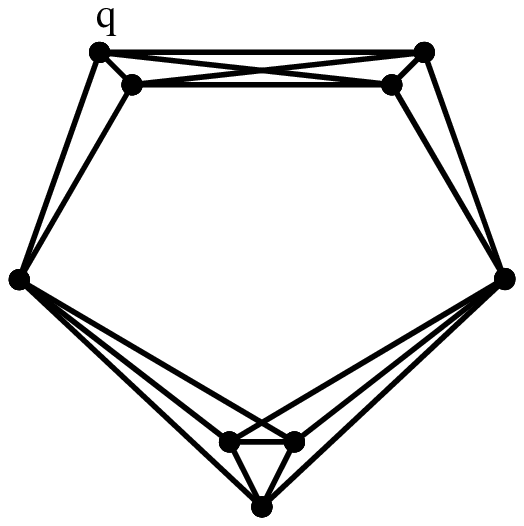}
\end{center}
\caption{}
\label{cutsetcountere}
\end{figure}

At first this result may seem not too helpful, because many of the graphs
for which the SSH conjecture is most mysterious have vertex degrees $
n-o(n)$ and connectivity $n-o(n)$. Nevertheless, it appears
this lemma does away with some pathological cases that would otherwise
present problems for a nice proof of the general result. In fact, we
conjecture that if $G$ has no cutset of size $\frac{n}{2}$ or smaller and
no dominating edge, then $G$ has a minor with any vertex $q\in V$ as
a prevertex and all the other prevertices as edges.  We can easily find graphs with no such minor, but all of them we found have dominating edges or a small cut set. For example, in the graph in figure \ref{cutsetcountere}, there is no $K_5$ minor that uses $q$ as a prevertex.

\subsection{Random graphs}

Let $G(n,p)$ be the Erd\"os-Renyi random
graph in which there are $n$ vertices and edges are independently present
with probability $p$. For every constant $p$, $0<p<1$, Hadwiger's conjecture is true for almost all graphs in $G(n,p)$ \cite{BCE}.  It is unlikely (though possibly still worth thinking about) that a reasonable random graph model will yield a counterexample to Hadwiger's conjecture or SSH.
Nonetheless, random graphs provide an aid to our intuition by helping us think
about graphs we cannot construct. We show some computations that suggest SSH holds for almost all graphs in $G(n,p)$ where $p = 1 - c n^{-\alpha}$, $\frac{1}{4} \leq \alpha$.

\begin{claim}
Let $p = 1 - c n^{-\alpha}$, $\frac{1}{4} \leq \alpha$, $q = 1-p$, and $d=\lfloor (n-1)/2 \rfloor$. For a graph in $G(n,p)$, the expected number of $K_{d+1}$ minors using only sets of size one or two as prevertices tends to infinity as $n$ tends to infinity.
\end{claim}

\begin{proof}
If $n$ is odd, Let $V(G) = \{v_1, v_2, \ldots, v_n\}$. The
probability that $\{v_1, v_2\}$, $\{v_3, v_4\}$, $\ldots$, $\{v_{n-2}, v_{n-1}\}$,
and $v_n$ are the prevertices of a complete minor is
\begin{equation*}
p^{d} (1-q^4)^{\binom{d}{2}} (1-q^2)^{d}
\end{equation*}
The terms in this product are the probability that $(v_{2i-1},v_{2i})$ is
an edge, that the prevertices of size 2 touch each other, and that $v_n$ touches the
prevertices of size 2.  Substituting $c n^{-\alpha}$ for $q$ we obtain
\begin{equation*}
(1 - c n^{-\alpha})^{d} (1-c^4 n^{-4\alpha})^{\binom{d}{2}} (1-c^2 n^{-2\alpha})^{d} . \end{equation*}
Using the expansion $\log(1-\epsilon) = -\epsilon - \text{$\frac{1}{2}$}\epsilon^2 - \ldots = -\epsilon - o(\epsilon)$ we obtain
$$\exp \left( -c n^{-\alpha}d - o(n^{1-\alpha}) - c^4 n^{-4\alpha} \binom{d}{2} - o(n^{2-4\alpha}) - c^2 n^{-2\alpha} d - o(n^{1-2\alpha}) \right) =  $$
\begin{equation}
\label{probeq}
\exp \left( - \frac{c}{2} \; n^{1-\alpha} - o(n^{1-\alpha}) - \frac{c^4}{8} \; n^{2-4\alpha} - o(n^{2-4\alpha}) \right) \geq e^{- \frac{c^4}{8} \; n + o(n)}
\end{equation}
since $2-4\alpha \leq 1$. There are $$ n \frac{1}{d!} \binom{n-1}{2,2, \ldots,2} = \frac{n!}{d! 2^d} $$
distinct sets of prevertices of the above type.  Using Stirling's approximation, there is a constant $c'$ so that this is
$$ \geq \sqrt{\frac{n}{d}} \frac{c' n^n}{e^{n-d} (2d)^d} = $$
\begin{equation}
\label{probeq2}
 \exp \left( n \log n - d \log 2d + d - n + \log (c' \sqrt{\frac{n}{d}}) \right) = e^{ \frac{n}{2} \log n + o(n \log n) }
\end{equation}
Combining (\ref{probeq}) and (\ref{probeq2}), the expected number of complete minors of this type is at least
$$ e^{ \frac{n}{2} \log n - \frac{c^4}{8} \; n + o(n \log n)} $$
which tends to infinity as $n$ goes to infinity.  A similar (and slightly simpler) argument works for $n$ even.
\end{proof}

Showing that SSH holds for almost all graphs (with $p$ as above) requires a second moment calculation.  This seems doable but tedious, and we do not do it.  We have yet to say anything about graphs with no antitriangle.  The expected number of antitriangles in $G \in G(n,p)$ is $\binom{n}{3} q^3 \sim \frac{c^3}{6} n^{3-3\alpha}$. For there to be asymptotically no antitriangles, we must have $\alpha > 1$.  A trick from Ramsey theory is to add edges to destroy antitriangles; this does not mess up the graph's properties too much if $\alpha \geq \frac{1}{2}$ \cite{Sudakov}. As $\alpha$ is decreased, such a strategy becomes less effective, and this graph model can say little about graphs with no antitriangle.  It therefore need not worry us that when $\alpha < \frac{1}{4}$, the expected number of minors of the type above tends to 0.

\subsection{Constant factor weakenings are unsolved}
One approach to Hadwiger's conjecture for graphs with no antitriangle is to try to show there is a complete minor of size $cn$ for some constant $c>0$, rather than demanding $c=1/2$. Even this weakening is unsolved for SSH.  We present the progress made in this direction, and begin with an instructive result observed independently by Mader, Kelmans, and Seymour.

\begin{claim}
If $G = (V,E) \in \mathfrak{A}$, then $G$ has a $K_{n/3}$ minor.
\end{claim}

\begin{proof}
We can obtain such a minor using induced paths of length 2 and
single vertices as prevertices. If $u, v, w$ are the vertices of an induced path of length 2, because there
is no antitriangle, $N(\{u,w\})$ = $V-\{u,w\}$. Choose a maximum
number of vertex disjoint induced paths of length 2. Let their vertex sets
be $Q_1, Q_2, \ldots Q_r$ and let $Q = \bigcup_i{Q_i}$. In $G\backslash Q$,
there are no induced paths of length 2, so being connected by an edge is an
equivalence relation. Thus $G\backslash Q$ is the disjoint union of at most
two cliques; let $C$ be the largest clique of $G\backslash Q$. $Q_1, Q_2,
\ldots Q_r$ and the vertices of $C$ are the prevertices of a $K_{r+|V(C)|}$
minor. $3r+2|V(C)| \geq n$ implies $r+|V(C)| \geq n/3$.
\end{proof}

Induced paths of length 2 are a bit of a cheat because they let us ignore the complex structure of these graphs.  For this reason the following problems are of interest.

\begin{problem*}
Show that there is a constant
\setcounter{ctr}{0}
\begin{list}{\emph{(\roman{ctr})}} {\usecounter{ctr} \setlength{\itemsep}{1pt} \setlength{\topsep}{2pt}}
\item $c>1/3$ such that for every $G \in \mathfrak{A}$, $G$ has a $K_{cn}$ minor.
\item $c>0$ such that for every $G \in \mathfrak{A}$, $G$ has a $K_{cn}$ minor using only cliques as prevertices.
\item $c>0$ such that for every $G \in \mathfrak{A}$, $G$ has a $K_{cn}$ minor using prevertices of size one or two.
\end{list}
\end{problem*}

Using an elementary counting argument, we show problem (iii) holds if $K_{cn}$ is replaced by $K_{cn^{4/5}}$.

\begin{theorem}
\label{n4/5claim} Let $G \in \mathfrak{A}$ have minimum degree $\delta(G) = n - c_1 n^{\alpha}$.  Assume that $0 \leq \alpha < 1$ so that $|E(G)| = \frac{1}{2} n^2 + o(n^2)$.  Then $G$ has a complete minor of size $c_3 n^{\beta} + o(n^{\beta})$ using prevertices of size one or two, where $\beta = \min(4-4\alpha,1)$ and $c_3$ is a constant depending only on $c_1$.
\end{theorem}

\begin{proof}
Let $H$ be the graph with vertex set $E(G)$; edges $e_1$ and $e_2$ are
adjacent in $H$ if they share an end or do not touch. A stable set in $H$
gives the prevertices of a complete minor in $G$. We will bound the degree
of $H$ to show that it has a large stable set.

If $e$ is an edge, let $\overline{N(e)}$ be the set of vertices that do not
touch $e$. A \emph{vedge} is the simple graph with three vertices and one edge. We
count the number of induced vedges in $G$ in two different ways.
\begin{equation}  \label{twocountings}
\sum_{v \in V(G)} \binom{n-d(v)}{2} = {\QATOP{\text{number of} }{\text{%
induced vedges}}} = \sum_{e \in E(G)} |\overline{N(e)}|
\end{equation}
$\binom{n-d(v)}{2}$ is the number of vedges with isolated vertex $v$, and $|%
\overline{N(e)}|$ is the number of vedges with edge $e$. Using the degree
bound, we obtain
\begin{equation}  \label{eq2}
n \; \frac{c_1^2}{2} \; n^{2 \alpha} \geq \sum_{v \in V(G)} \binom{n-d(v)}{2}
\end{equation}

Then the average value of $|\overline{N(e)}|$ is about $c_1^2 n^{2\alpha-1}$.  Let $E^{\prime}$ be the edges $e$ for which $|\overline{N(e)}| \geq 2 c_1^2 n^{2\alpha-1}$ (twice the average is arbitrary; other constant factors would do). We may now bound $|E^{\prime}|$.  Define $c_2$ so that $|E^{\prime}|= c_2 |E(G)|$. Then by (\ref{twocountings}) and (\ref{eq2})
\begin{equation*}
\frac{c_1^2}{2} \; n^{2 \alpha+1} \geq \sum_{e \in E(G)} |\overline{N(e)}%
| \geq c_2 |E(G)| \; 2 c_1^2 \; n^{2\alpha-1}
\end{equation*}
implies
\begin{equation} \label{eq3}
c_2 \leq \frac{n^2}{4 |E(G)|} = \frac{1}{2} + o(1)
\end{equation}
Then for $|E(G)-E'| \geq (1-(\frac{1}{2}+o(1)))|E(G)|$ edges $e$,
\begin{equation*}
d_H(e) \leq 2n + \binom {|\overline{N(e)}|}{2} \leq 2n + 2c_1^4
n^{4\alpha-2}.
\end{equation*}%
The bound on $d_H(e)$ comes from the trivial upper bound of $2n$ for the
number of edges sharing an end with $e$, and $\binom {|\overline{N(e)}|}{2}$
is from the fact that $\overline{N(e)}$ is a clique containing all edges not
touching $e$. Then $H\backslash E^{\prime}$ has max degree $\Delta \equiv 2n
+ 2c_1^4 n^{4\alpha-2}$ and a greedy coloring shows the chromatic number $%
\chi(H\backslash E^{\prime}) \leq \Delta+1$. This together with (\ref{eq3}) implies there is a stable
set in $H \backslash E^{\prime}$ of size at least $\frac{|E(G)|-|E^{\prime}|}{\Delta+1} \geq \frac{
(1/2 + o(1))|E(G)|}{\Delta + 1} = \frac{n^2/4}{\Delta}(1+o(1)) $

Put
\begin{equation*}
c_3=\left\{
\begin{array}{ll}
\frac{1}{8c_1^4} & \text{if } 4\alpha-2 > 1 \\
\frac{1}{4(2+2c_1^4)} & \text{if } 4\alpha-2 = 1 \\
\frac{1}{8} & \text{if } 4\alpha-2 < 1
\end{array}
\right.
\end{equation*}

Put $\beta = \min(4-4\alpha,1)$. Then $G$ has a complete minor of size $\frac{n^2/4}{\Delta}(1+o(1)) = c_3
n^{\beta} + o(n^{\beta})$.
\end{proof}

The constants obtained in the proof are not the optimal obtainable by this method, but they will do.  The corollary below follows easily.
\begin{corollary}
For sufficiently large $n$, every $G \in \mathfrak{A}$ has a complete minor of size $\frac{n^{4/5}}{9}$ using prevertices of size one or two.
\end{corollary}
In a graph with no antitriangle the non-neighbors of each vertex are a clique. Then $\delta(G) = n - c_1 n^{\alpha}$ implies $G$ has a complete minor of size $c_1 n^{\alpha}$. Note
that $\max(\min(4-4\alpha,1),\alpha) \geq \frac{4}{5}$.  Also observe that $\max(\frac{1}{4(2+2c_1^4)},c_1) > 1/9$ and the corollary follows.

This method shows $G$ has a complete minor of size $O(n)$ when $\alpha \leq \frac{3}{4}$, that is, when $\delta(G) \geq n - c_1 n^{3/4}$.  Random graphs, we suspect, have a complete minor when $p \geq 1 - c n^{-1/4}$, that is, when the expected degree of a vertex is $\geq n - c n^{3/4}$.  That this threshold is the same is interesting, but probably coincidental.

\section{Good and bad edges}

Let $c_1, c_2, \ldots c_r$ be the cliques of $G$ and let $w$ be a function from $\{c_i\}$
to the nonnegative rationals. The \emph{fractional clique covering number} of $G$
is the minimum of $\sum_i w(c_i)$ over all maps $w$ such that
\begin{equation*}
\forall \; v \in V(G) \; \sum_{i \; s.t. \; v \in \; \large{c_i}} w(c_i) \geq 1 .
\end{equation*}
If $G$ has fractional clique covering number less than 3, multiplying $w$
by a common denominator shows that there is a list of $k$ cliques (not necessarily distinct) such that every vertex is in more than $\frac{k}{3}$ of them.
In particular this implies that $G$ has no antitriangle.
It is interesting to study the SSH conjecture for such graphs.

We observe that there is a natural way to partition the edges in a graph with fractional clique covering number less than 3.
An edge $(u,v)$ is \emph{good} if there are more than $\frac{k}{2}$ cliques
containing $u$ or $v$. If $(u,v)$ and $(x,y)$ are good, there is a clique
containing at least one of $u,v$ and at least one of $x,y$, so every pair of
good edges touch. An edge $(u,v)$ is \emph{bad} if there are $\frac{k}{2}$
or fewer cliques containing $u$ or $v$. If $(u,v)$ and $(v,w)$ are bad, then
there are $<\frac{k}{6}$ cliques containing $v$ and not $u$ and $<\frac{k}{6}
$ cliques containing $v$ and not $w$, so there is a clique containing $%
\{u,v,w\}$. In other words, there are no induced paths of length 2 that use only bad edges.
It turns out that under some not too restrictive conditions, $G$ has a
perfect matching of good edges, and these edges are the prevertices of a complete
minor. All that is needed to prove this are the conditions on pairs of
good and bad edges. We therefore drop the fractional clique covering number
condition and retain the conditions on edge pairs.

\subsection{Perfect matching of good edges}

Let $G=(V,E)$ be a graph with no antitriangle. Suppose $E$ can be
partitioned into good edges and bad edges, $E=\mathfrak{G}\cup \mathfrak{B}$%
, so that for every pair $g_{1},g_{2}$ of good edges, $g_{1}$ and $g_{2}$
touch, and for every pair of bad edges that share an end, $b_{1}=\{u,v\}, \;
b_{2}=\{v,w\}$, $\{u,w\}$ is an edge. We will call these conditions the good
edge and bad edge axioms.

\begin{theorem}
\label{Goodbad} If $G = (V, E) \in \mathfrak{A}$, $n$ is even, and $E=\mathfrak{G}\cup
\mathfrak{B}$ as above, then $G$ has a $K_{n/2}$ minor with prevertices of
size at most 2.
\end{theorem}

This subsection and part of the next are devoted to proving this theorem. If $G$ has
a dominating edge $(u,v)$, then $E(G \backslash \{u,v\})$ can also
be partitioned into good edges and bad edges. By induction on $n$, we
obtain a complete minor of $G \backslash \{u,v\}$; by adding the prevertex $%
(u,v)$, we obtain a complete minor of $G$. We know the minor exists if there
is a small cutset (lemma \ref{cutset}). Obviously, the minor exists if $G$
has a clique of size $\geq n/2$. If none of these arguments works, we
prove that $G$ has a perfect matching of good edges unless $G$ is the
complement of a blown up Petersen graph. This gives the prevertices of a $%
K_{n/2}$ minor. If $G$ is the complement of a blown up Petersen graph, the
minor can be found easily in several ways.

\begin{theorem}
\label{goodmatching} If $G = (V, E) \in \mathfrak{A}$, $n$ is even, and $E=\mathfrak{G}\cup
\mathfrak{B}$ as above, then either $G$
\setcounter{ctr}{0}
\begin{list}{\emph{(\alph{ctr})}} {\usecounter{ctr} \setlength{\itemsep}{1pt} \setlength{\topsep}{2pt}}
\item has a dominating edge, or
\item has a cut set of size $\leq n/2$, or
\item has a clique of size $\geq n/2$, or
\item is the complement of a blown up Petersen graph, or
\item has a perfect matching of good edges.
\end{list}
\end{theorem}

\begin{proof}
We will assume (a), (b), (c), and (e) are false and prove (d). Let $G^{\prime}
$ $\equiv (V,\mathfrak{G})$ and let $\mathfrak{N}=E(\overline{G})$. We apply
Tutte's theorem to $G^{\prime}$: (e) is false implies $\exists $ $S\subseteq
V$ such that $G^{\prime }\backslash S$ has at least $|S|+2$ odd components
(components with an odd number of vertices). We will mostly be working with
the graph $G$, and often think of it as a complete graph with three edge
types--good ($\mathfrak{G}$), bad ($\mathfrak{B}$), and non-edges ($%
\mathfrak{N}$). We will occasionally refer to $G^{\prime}$; be sure not to
confuse the two.

Let $C_{1},C_{2},\ldots ,C_{m}$ be the components of $G^{\prime }\backslash
S $; $m \geq |S|+2$.

\begin{step}{}
 Either
\setcounter{ctr}{0}
\begin{list}{\emph{(\roman{ctr})}} {\usecounter{ctr} \setlength{\itemsep}{1pt} \setlength{\topsep}{2pt}}
\item $\overline{G} \backslash S$ is bipartite, or
\item $C_i$ contains an odd antihole, for some $i \in [m]$, or
\item $G$ contains an antihole of length 5 that belongs to two components and is
isomorphic to the 5-antihole in figure \ref{5cycle}.
\end{list}
\end{step}

\par
If $\overline{G}\backslash S$ is not bipartite, it contains an odd cycle of
length 5 or greater (a cycle of length 3 is an antitriangle in $G$).
A shortest such cycle is an odd antihole in $G\backslash S$, which we label $%
Y=v_{1},v_{2},\ldots ,v_{l}$. Throughout the proof of (1), we treat all
subscripts $mod \; l$. Let
\begin{equation*}
D_{i,j}=
\left\{ \begin{array}{ll}
1 & \text{if $v_i$ and $v_j$ belong to the same component} \\
0 & \text{otherwise}
\end{array} \right.
\end{equation*}
We will call $v_{i},v_{j},v_{j+1}\in Y$, a \emph{%
forcing triple} if $d(i,$ $j)$, $d(i,$ $j+1)>1$, where $d$ is the distance $
mod \; l$. Edges $(v_{i},v_{j})$ and $(v_{i},v_{j+1})$ are not both bad, so $%
D_{i,j} = 1$ or $D_{i,j+1}=1$. Either $l \geq 7$ (case A), or $l = 5$ (case
B).

\medskip \noindent (A) First suppose $Y$ intersects at most two
components. Then $\exists \; j$ such that $D_{j,j+1} = 1$, and therefore $%
D_{i,j} = 1$ for all $i \in [l]-\{j-1,j+2\}$. The forcing triples $%
v_{j+3},v_{j+4},v_{j-1}$ and $v_{j-2},v_{j-3},v_{j+2}$ show that $D_{j,j-1}
= D_{j,j+2} = 1$. Thus $Y$ is contained one component ((ii) holds). If $Y$
intersects more than two components, just merge all but one of them and
treat it as a single component. The same proof works.

\medskip \noindent (B) Begin as in (A) by supposing $Y$ intersects at most
two components, and choose $j$ as in (A). The forcing triple $%
v_{j-1},v_{j},v_{j+2}$ shows $D_{j,j+2}=1$ or $D_{j-1,j+2}=1$. If the
former, the forcing triple $v_{j+1},v_{j+2},v_{j-1}$ shows (ii) holds. If
the latter holds, but the former does not, then (iii) holds. If $Y$
intersects more than two components, then without loss of generality, $C_1
\cap Y = v_1$. The forcing triple $v_1,v_3,v_4$ gives a contradiction.

\begin{figure}[tbp]
\begin{center}
\includegraphics [scale=1]{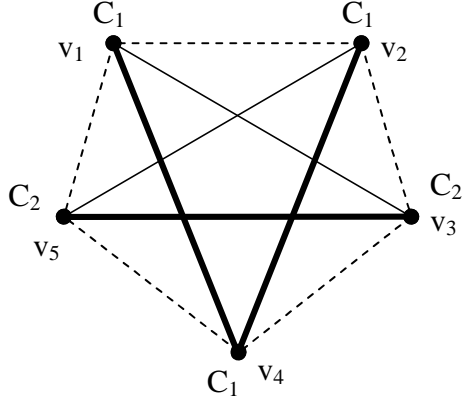}
\end{center}
\caption{The antihole of (1)(iii). The thicker (thinner) edges are good
(bad), and the dotted edges are non-edges. The labeling of good and bad
edges is forced by the component assignments.}
\label{5cycle}
\end{figure}

\bigskip

The longest part of the proof is the $m = 2$ case, which we treat specially.
Steps (2) and (3) are devoted to this case and step (4) addresses the $m > 2$ case.
For convenience, let $L=C_1$, $R=C_2$. Before proceeding, we need some definitions.

A set of vertices $P$ is $X$-\emph{coupled}, $X\subseteq V$, if $\forall \; p_1, p_2 \in P$, $N(p_1)\cap X=N(p_2)\cap X$. If $P$ and $Q$ are vertex
sets that are both $X$-coupled, we say $P,Q$ is $X$-\emph{anticoupled}, if $%
N(P)\cap X = X-N(Q)$. Also, we say an edge $(u,v)$ is $X$-coupled if $\{u,
\; v\}$ is $X$-coupled and $X$-anticoupled if $u,v$ is $X$-anticoupled.

Let $M_{1}$ be a component of bad and non-edges in $L$ (a component in the
graph $(V,\mathfrak{B} \cup \mathfrak{N}$)). Only bad edges and non-edges
cross between $L$ and $R$; by the bad edge axiom and no antitriangle,
every bad edge in $M_{1}$ is $R$-coupled and every non-edge in $M_{1}$ is $R
$-anticoupled. Thus there is no odd-cycle of non-edges, and $M_{1}$ can be
partitioned into two sets, $M_{1T}$ and $M_{1B}$, so that $M_{1T}, M_{1B}$
is $R$-anticoupled. Note that $M_{1T}$ and $M_{1B}$ are cliques (in $G$). We
will call $M_1$ a \emph{dipole} and call $M_{1T}$ and $M_{1B}$ \emph{poles}.
Given two poles of a dipole, we say one is the \emph{antipole} of the other.
If both (exactly one) poles of a dipole are nonempty, we will say the dipole
is \emph{proper} (\emph{improper}).

Let $l$ ($r$) be the number of dipoles in $L$ ($R$). We have $L=M_{1}\cup
M_{2}\ldots \cup M_{l}$ and $R=N_{1}\cup N_{2},\ldots \cup N_{r}$. By
definition of the $M_{i}$, any edge $\{u,v\}$ with $u\in
M_{i},v\in M_{j}$ $(i\neq j)$ is good. An edge between $M_{iT}$ and $%
M_{iB}$ is dominating because it touches all of $L-M_{i}$ by the good edges
just mentioned and touches all of $R$ because it is $R$-anticoupled. Thus a
pole does not touch its antipole (in $G$) and every pair of vertices in a pole are c-twins.
For every $i \in [l], \; j
\in [r]$, $M_{iT}$ touches exactly one of $N_{jT}, N_{jB}$, because $%
N_{jT}, N_{jB}$ are $L$-anticoupled. So either ($M_{iT}$ is complete to $%
N_{jT}$ and $M_{iB}$ is complete to $N_{jB}$) or ($M_{iT}$ is complete to $%
N_{jB}$ and $M_{iB}$ is complete to $N_{jT}$)---if the former, we say the
dipoles are \emph{matched straight} and if the latter they are \emph{%
matched twisted}.

If 1(i) holds, there is a large clique, but we assumed (c) is false,
contradiction. The decomposition into dipoles shows that $\overline{G}%
\backslash L$ and $\overline{G}\backslash R$ are bipartite so 1(ii) does not
hold. We may assume 1(iii) holds. From figure \ref{5cycle}, we see that $%
v_3, v_5$ are neither $L$-coupled nor $L$-anticoupled, so $r \geq 2$. $v_1,
v_2$ are $R$-anticoupled, and $v_1, v_4$ are neither $R$-coupled nor $R$%
-anticoupled. So $l \geq 2$ and at least one dipole in $L$ is proper.

\begin{step} 
If two dipoles of $L$ are proper, then either
\setcounter{ctr}{0}
\begin{list}{\emph{(\roman{ctr})}} {\usecounter{ctr} \setlength{\itemsep}{1pt} \setlength{\topsep}{2pt}}
\item $G$ is the complement of a blown up $V_8$ as shown in figure \ref{v8complement}, or
\item $G$ is the complement of a blown up Petersen graph as shown in figure \ref{Petersencomplement}.
\end{list}
\end{step}

Without loss of generality, $M_1$ and $M_2$ are proper. Consider \newline $%
N(M_{1T})\cap R, \; N(M_{1B})\cap R, \; N(M_{2T})\cap R, \; N(M_{2B})\cap R$%
, and call them $T_1,B_1,T_2,B_2$ for brevity. Every pole in $R$ is in
exactly two of these sets. An edge, $e$, between $M_{1T}$ and $M_{2T}$ is
good and therefore touches every good edge, so at most one dipole intersects
$R - (T_1 \cup T_2) = B_1 \cap B_2$. If no dipole intersects $B_1 \cap B_2$,
then $e$ is dominating, so we may assume exactly one dipole intersects $B_1
\cap B_2$. Applying the same argument to edges between $M_{1T}$ and $M_{2B}$%
, $M_{1B}$ and $M_{2T}$, and $M_{1B}$ and $M_{2B}$ shows that exactly one
dipole intersects $B_1 \cap T_2$, $T_1 \cap B_2$, and $T_1 \cap T_2$.

If $R$ contains a proper dipole, $N_1$, say, then it must have non-empty
intersection with each of $T_1,B_1,T_2$, and $B_2$ (because $N_{1T}, N_{1B}$
is $L$-anticoupled). Then (up to symmetry between $N_{1T}$ and $N_{1B}$)
either ($N_{1T} = T_1 \cap T_2$ and $N_{1B} = B_1 \cap B_2$) or ($N_{1T} =
T_1 \cap B_2$ and $N_{1B} = T_2 \cap B_1$). Clearly, since every vertex of a
pole in $R$ has the same neighbors in $L$, poles are contained in the sets $%
B_1 \cap B_2$, etc. As just seen, both poles of a dipole cannot be contained
in $B_1 \cap B_2$, etc., so, in fact, $B_1 \cap B_2, B_1 \cap T_2$, $T_1
\cap B_1$, and $T_1 \cap T_2$ are poles. Up to symmetry, there are three
possibilities, (A), (B) and (C), for the structure of $R$.

\medskip \noindent (A) $R$ is the union of two proper dipoles: $N_{1T} =
T_1 \cap T_2$, $N_{1B} = B_1 \cap B_2$, $N_{2T} = T_1 \cap B_2$, and $N_{2B}
= T_2 \cap B_1$. Applying the argument above with $L$ and $R$ reversed,
shows that $L$ is the union of four poles, which implies $l=2$. We now know the structure of $G$ up to vertex c-duplication---$G$ is the complement of a blown
up $V_8$ ((i) holds).

\begin{figure}[tbp]
\begin{center}
\includegraphics [scale=1.2]{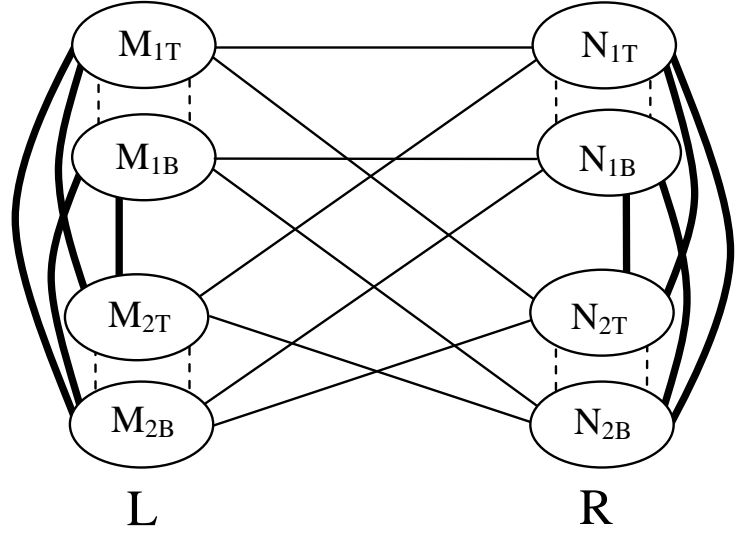}
\end{center}
\caption{The complement of a blown up $V_8$. All good and bad edges are drawn. Non-edges between poles and antipoles are drawn, but non-edges between $L$ and $R$ are not.}
\label{v8complement}
\end{figure}

\medskip \noindent (B) $R$ is the union of a proper dipole and two
improper dipoles: $N_{1T} = T_1 \cap T_2$, $N_{1B} = B_1 \cap B_2$, $N_{2} =
T_1 \cap B_2$, and $N_{3} = T_2 \cap B_1$. An $(N_2, N_3)$ edge is not
dominating, so $l \geq 3$. Furthermore, there is a pole, $M_{3T}$, say, that
doesn't touch $N_2$ or $N_3$. Without loss of generality, $M_{3T}$ touches $%
N_{1T}$ and not $N_{1B}$. But then an $(N_2,N_{1B})$ edge doesn't touch an
$(M_{2T}, M_{3T})$ edge, contradicting the good edge axiom.

\begin{figure}[tbp]
\begin{center}
\includegraphics [scale=1]{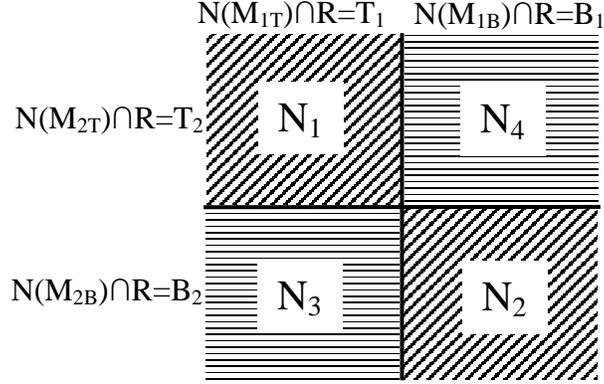}
\end{center}
\caption{Schematic drawing of the partitions $T_1 \uplus B_1$ and $T_2 \uplus B_2$.  The different shadings represent the third partition $N(M_{3T})\cap R \uplus N(M_{3B})\cap R$.}
\label{partition}
\end{figure}

\medskip \noindent (C) $R$ is the union of four improper dipoles: $N_{1} =
T_1 \cap T_2$, $N_{2} = B_1 \cap B_2$, $N_{3} = T_1 \cap B_2$, and $N_{4} =
T_2 \cap B_1$ (see figure \ref{partition}). An $(N_1, N_2)$ edge is not dominating so $l \geq 3$. Any
dipole in $L$ is proper, because suppose $M_3$ were improper. Then an edge
from $M_3$ to $R$ is dominating; if $M_3$ doesn't touch $R$, the good edge
axiom is violated. Now apply the same argument to $\{M_1, M_j\}$ and $\{M_2,
M_j\}$ as was applied to $\{M_1, M_2\}$, $j \geq 3$. Each dipole in $L$
partitions $R$ into two sets each containing two poles. Moreover, every two
such partitions must be isomorphic to the partitions defined by $M_1$ and $%
M_2$ ($T_1 \uplus B_1$ and $T_2 \uplus B_2$). That leaves room for only one
more partition: $(N_1 \cup N_2) \uplus (N_3 \cup N_4)$. So $l=3$ and without
loss of generality, $N(M_{3T})\cap R = N_1 \cup N_2$ and $N(M_{3B})\cap R = N_3 \cup N_4$%
. We now know the structure of $G$ up to vertex c-duplication---$G$ is the
complement of a blown up Petersen graph ((ii) holds). This proves (2).

\begin{figure}[tbp]
\begin{center}
\includegraphics [scale=1]{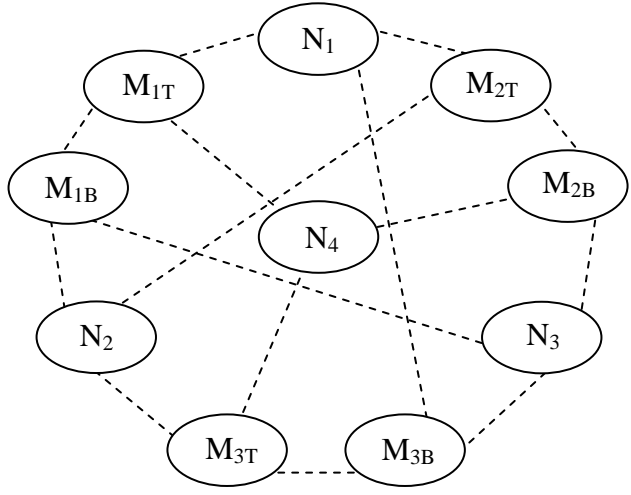}
\end{center}
\caption{The complement of a blown up Petersen graph.  All non-edges are drawn.}
\label{Petersencomplement}
\end{figure}

If $R$ has two proper dipoles, we may apply (2) and conclude (e). So we may
assume $L$ and $R$ have at most one proper dipole. We already know $L$ has
at least one proper dipole, so $L$ has exactly one proper dipole. We now
know a lot about the structure of $G$, and can finish up the remaining cases in
(3).

\begin{step}
Given the conclusions of (1) and (2), we may assume
\setcounter{ctr}{0}
\begin{list}{\emph{(\roman{ctr})}} {\usecounter{ctr} \setlength{\itemsep}{1pt} \setlength{\topsep}{2pt}}
\item $L$ contains exactly one proper dipole,
\item $R$ contains at most one proper dipole,
\item $l, r \geq 2$,
\item $L$ is the union of a proper dipole and an improper dipole, and
\item $G$ is isomorphic to a graph represented by figure \ref{v8complement} with $M_{2B}=N_{2B}= \emptyset$ and all other sets nonempty except possibly $N_{1T}$.
\end{list}
\end{step}

As just discussed, (i) and (ii) hold. (iii) we have already seen. Without
loss of generality, $M_1$ is proper. We proceed as in the proof of (2).
There is less symmetry so the arguments are a bit messier. Consider $%
N(M_{1T})\cap R, \; N(M_{1B})\cap R, \; N(M_2)\cap R, \; R-N(M_2)$, and call
them $T_1,B_1,T_2, B_2$ for brevity. An edge, $e$, between $M_{1T}$ and $M_2$
is good and therefore touches every good edge, so at most one dipole
intersects $R - (T_1 \cup T_2) = B_1 \cap B_2$. If no dipole intersects $B_1
\cap B_2$, then $e$ is dominating, so exactly one dipole intersects $B_1
\cap B_2$. Applying the same argument to edges between $M_{1B}$ and $M_2$
shows that exactly one dipole intersects $T_1 \cap B_2$. Either $R$ contains
two improper dipoles (case A), or it does not (case B).

\medskip \noindent (A) $R$ contains at least two improper dipoles, $N_1$
and $N_2$. $N_1$ is not contained in $T_2$ because then an $(M_2,N_1)$ edge
is dominating. Similarly for $N_2$. From the discussion above, we must have
(up to symmetry of labeling) $N_{1} = T_1 \cap B_2$, and $N_{2} = B_1 \cap
B_2$. Suppose for a contradiction that $l \geq 3$. $M_3$ is improper so any
edge from $\{M_2 \cup M_3\}$ to $\{N_1 \cup N_2\}$ is dominating. If $\{M_2
\cup M_3\}$ does not touch $\{N_1 \cup N_2\}$ this violates the good edge
axiom. So (iv) holds. If $R$ contains another dipole, $N_3$, it is contained
in $T_2$, but then either an $(N_1,N_3)$ or an $(N_2,N_3)$ edge is
dominating, contradiction. We have determined the structure of $G$ up to vertex
c-duplication---$G$ is isomorphic to figure \ref{v8complement} with $%
M_{2B}=N_{2B}=N_{1T}=\emptyset$ and all other sets nonempty ((v) holds).

\medskip \noindent (B) By (ii) and (iii), $R$ is the union of a proper
dipole and an improper dipole. Apply (3) with $L$ and $R$ reversed. If $L$
contains more than one improper dipole, this is dealt with by (A). So we may
assume $L$ has exactly one proper dipole; this together with (i) implies
(iv). The proper dipole of $R$, $N_1$, say, must have non-empty intersection with
each of $T_1,B_1$ and $T_2$ (because $N_{1T}, N_{1B}$ is $L$-anticoupled).
We know from discussion above that $T_1 \cap B_2$ and $B_1 \cap B_2$ are
poles. This determines the structure of $G$ up to vertex c-duplication---$G$ is
isomorphic to figure \ref{v8complement} with $M_{2B}=N_{2B}=\emptyset$ and all
other sets nonempty ((v) holds). This proves (3). 

\medskip
If (2)(i) or (3)(v) holds, $G\backslash (M_1 \cup N_1)$ and $G\backslash
(M_2 \cup N_2)$ are disconnected. At least one of $|M_1 \cup N_1|$ and $|M_2
\cup N_2|$ is $\leq n/2$ so (c) is true, contradiction. If (2)(ii) holds,
(d) is true, as desired. This completes the $m=2$ case.

\bigskip
We may assume $m > 2$. It is here we reap the main rewards of (1). If
(1)(ii) or (1)(iii) holds, then an antihole, $Y$, does not intersect all
components of $G^{\prime} \backslash S$. Let $L$ ($R$) be the union of all
components $Y$ intersects (doesn't intersect). Non-edges in $A$ are $B$%
-anticoupled and therefore, as seen earlier, an odd-cycle of non-edges is
impossible. So we may assume (1)(i).

(1)(i) implies we can partition $V-S$ into $A$ and $B$ such that $
(G\backslash S)[A]$ and $(G\backslash S)[B]$ are cliques. Let $C_{i}^{A}=C_{i}\cap A$ and $C_{i}^{B}=C_{i}\cap B$.
Since $|C_{i}|$ is odd, $|C_{i}^{A}|\neq |C_{i}^{B}|$ . Let $X=$ $\bigcup _{i}$%
(smaller of $C_{i}^{A}$, $C_{i}^{B}$). Remembering that $m \geq |S| +2$, we
observe $|X\cup S| \leq n/2$. Therefore $X\cup S$ is not a cutset and $%
G\backslash X\backslash S$ is connected.

Without loss of generality $G\backslash X\backslash S$ is the union of $%
C_{1}^{A},C_{2}^{A},...$, $C_{k}^{A}$, $C_{k+1}^{B}$, $C_{k+2}^{B},...,C_{m}^{B}$,
$0 \leq k \leq m$. By symmetry we may assume $k \leq m-k$.

\begin{step}
The cases $k=0$ (A), $k \geq 2$ (B), and $k=1$ (C) each lead to a contradiction, which shows that $m > 2$ is impossible.
\end{step}

\noindent (A) $G\backslash X\backslash S$ is a clique and it is large enough
to contradict the assumption that (c) is false.

\medskip \noindent (B) We may view $A-X$ and $B-X$ as $L$ and $R$ in the $%
m=2$ case because no good edges have one end in $A-X$ and one end in $B-X$
(remember, the $C_{i}$ are components of good edges). Since $k\geq 2$, all
of $A-X$ is a component of bad edges ($A-X$ is a pole). Since $m-k \geq 2
$, $B-X$ is a pole. $G\backslash X\backslash S$ is connected, so there is an
edge between $A-X$ and $B-X$, and therefore $A-X$ and $B-X$ are joined
completely. This proves (c), which we assumed false.

\medskip \noindent (C) $m > 2$ implies $m-k\geq 2$. For this case, we
will apply dipole structure to the partition $B-X, C_1$ (each $(B-X,C_1)$
edge is bad). As in (B), $B-X$ is a pole. Let $M_{B}=C_{1}^{A}\cap N(B-X)$. $%
C_{1}^{B}$ is complete to $M_{B}$ (in $G$) by the bad edge axiom. Since $|C_1|$ is odd, $|C_1^A-M_B|$ and $|C_1^B \cup M_B|$ are not equal.  If $|C_1^A-M_B|$ is larger, then $G \backslash (X \cup M_B \cup S)$ is disconnected and this contradicts the assumption that (b) is false; If $|C_1^B \cup M_B|$ is larger, then $B \cup M_B$ is a clique and this contradicts the assumption that (c) is false.  This proves (4).
\end{proof}

\subsection{Extensions}

We first give another, quite simpler, proof of theorem \ref{Goodbad} by
modifying theorem \ref{goodmatching}.

\begin{corollary}
\label{goodmatchingcor} If $G = (V, E) \in \mathfrak{A}$, $n$ is even, and $E=\mathfrak{G}\cup
\mathfrak{B}$ such that
\setcounter{ctr}{0}
\begin{list}{\emph{(\alph{ctr})}} {\usecounter{ctr} \setlength{\itemsep}{1pt} \setlength{\topsep}{3pt}}
\item $\mathfrak{G}$ is chosen as large as possible with the restriction that
\item all edges between c-twins are bad,
\end{list}
then at least one of theorem \ref{goodmatching}(a), (b), (c), (e) holds.
\end{corollary}

\begin{proof}
Note that given any partition of $E$ into good and bad edges with the axioms
satisfied, an edge whose ends are c-twins can be made bad without violating
the axioms. Therefore, if there is some partition satisfying the axioms,
then there is one satisfying the axioms and (a) and (b).

We replace (2) and (3) by the following argument and leave the rest of the
proof the same. Because vertices in a pole are c-twins, all edges with both
ends in the same pole are bad. We may assume $G$ contains a subgraph
isomorphic to the graph in figure \ref{5cycle}. $v_1$ and $v_3$ are not c-twins.
Let $M_{1T}$ be the pole containing $v_1$ and let $N_{1T}$ be the pole
containing $v_3$. Making edge $(v_1,v_3)$ good does not violate the good
edge axiom because all edges with both ends in $V - N(\{v_1,v_3\}) = M_{1B}
\cup N_{1B} $ are bad. This contradicts (a).
\end{proof}

This proof takes care of the case when $G$ is the complement of a blown up
Petersen graph. Note that the partition of edges for the complement of a blown up
Petersen graph that satisfies corollary \ref{goodmatchingcor}(a) and (b) is: all edges with both ends in a c-blown up vertex are bad and all other edges are good.

Let $G=(V,E)$ be a graph in $\mathfrak{A}$ and suppose $E$ is
partitioned, $E=\mathfrak{M}\cup \mathfrak{B}$, so that the bad edge
axiom holds for $\mathfrak{B}$, but the good edge axiom does not (necessarily) hold for $\mathfrak{M}$ (call them medium edges). Is it true that there is a perfect matching of medium edges? The complement of the Petersen graph with suitable vertex c-duplication is a counterexample to this question, but are there others?  While a
matching of medium edges would not necessarily give the prevertices of a complete minor, it
would be necessary for there to be a complete minor that does not use bad
edges as prevertices.
It might be useful to know when we can ignore some edges (edges that aren't adjacent to many edges, perhaps) and still find a perfect matching in the remaining edges.

\begin{figure}[tbp]
\begin{center}
\includegraphics [scale=1]{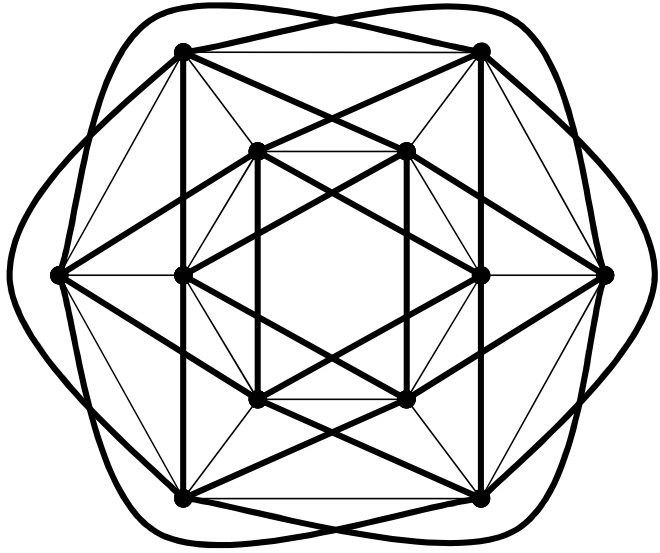}
\end{center}
\caption{The thicker (thinner) edges are medium (bad) edges.}

\label{12vcountere}
\end{figure}

There is not always a perfect matching of medium edges---we will see that a c-blown up version of figure \ref{12vcountere} is a counterexample. However, we have quite a bit of control on the counterexamples.  Note that the proof of theorem \ref{goodmatching} only uses the good edge axiom in steps (2) and (3).  So the only counterexamples
have the dipole structure described in the proof.  Moreover, $l, r \geq 2$, and $L$ contains at least one proper dipole.  Let $T$ be the complete bipartite graph with vertex set the set of dipoles.
If two dipoles are matched straight (twisted), label the corresponding edge
in $T$ straight (twisted). By exchanging labels of a pole and antipole, we
may swap the edge type of all edges incident to a vertex of $T$. Note that $%
T $, together with the number of vertices in every pole, is enough
information to reconstruct $G$, and graph(s) $T$ that yield a fixed $G$ are
not necessarily unique.  The graph in figure \ref{12vcountere} corresponds to the bipartite graph $T=K_{3,3}$ in which 3 vertex disjoint edges are twisted.  All poles are of size 1.  The smallest cutset in this graph has size 7.  This graph does have a perfect matching of medium edges, however, the graph with one pole in $L$ of size $k+1$ and one pole in $R$ of size $k+1$, and all other poles of size $k$ is a counterexample for large $k$.  This is because $L$ and $R$ are odd components in $(V,\mathfrak{M})$ so there is no matching of medium edges; the smallest cutset has size $7k$ which is $> n/2 = (12k+2)/2 $ for $k \geq 2$.

If we are willing to choose medium edges and bad edges with some additional properties, we can obtain a perfect matching of medium edges.  The trick from corollary \ref{goodmatchingcor} works with a simple modification.
\begin{corollary}
\label{goodmatchingcor2} If $G = (V, E) \in \mathfrak{A}$, $n$ is even, and $E=\mathfrak{M}\cup
\mathfrak{B}$ such that
\setcounter{ctr}{0}
\begin{list}{\emph{(\alph{ctr})}} {\usecounter{ctr} \setlength{\itemsep}{1pt} \setlength{\topsep}{3pt}}
\item all edges between c-twins are bad, and
\item given (a), the number of pairs of edges in $\mathfrak{M}$ that do not touch is as small as possible, and
\item given (a) and (b), $\mathfrak{M}$ is as large as possible,
\end{list}
then at least one of theorem \ref{goodmatching}(a), (b), (c), (e) holds (replace good with medium in (e)).
\end{corollary}

\begin{proof}
The proof is nearly the same as that of corollary \ref{goodmatchingcor}: we may assume $G$ contains a subgraph
isomorphic to the graph in figure \ref{5cycle}. $v_1$ and $v_3$ are not c-twins.
Let $M_{1T}$ be the pole containing $v_1$ and let $N_{1T}$ be the pole
containing $v_3$. Making edge $(v_1,v_3)$ good does not increase the number of pairs in $\mathfrak{M}$ that do not touch because all edges with both ends in $V - N(\{v_1,v_3\}) = M_{1B}
\cup N_{1B}$ are bad. Then $|\mathfrak{M}|$ was not maximum, contradicting (c).
\end{proof}

This corollary generalizes corollary \ref{goodmatchingcor} because if there is a way to partition the edges into good edges and bad edges, then the $\mathfrak{M}$ that satisfies (a), (b), and (c), will satisfy the good edge axiom.  Another way to look at the type of partition we are getting is that (roughly) the most useful edge sets are those where the good edge axiom is satisfied, so make such a set as large as possible.  Then, of the remaining, take an edge set that satisfies the bad edge axiom; it should be big so that the leftover edges, which must be added to $\mathfrak{M}$, do not ruin the good edge axiom too much.

Although the good edge axiom is what allows us to say anything about SSH, there is a reason we are trying to get rid of it in these generalizations.  It is too difficult to satisfy. It may be that in certain classes of graphs, finding a reasonable edge set satisfying the good edge axiom is as difficult as finding a complete minor. It is easier to identify
bad edges. If an edge, $e$, is between two vertices that are c-twins, then our
investigations strongly suggest we should be able to obtain a complete minor
without using $e$ as a prevertex.  Edges that connect two vertices that are ``close'' to being c-twins should also be labelled bad, with higher priority given to those that are closer. This will be discussed further in the conclusions.

\section{2 satisfiability}

We are fortunate that the good edge and bad edge axioms have a nice converse.
Every pair of edges that do not touch corresponds to a clause requiring that at
least one edge of the pair is bad. Every pair of edges $(u,v)$, $(v,w)$ such that $u \nsim w$
corresponds to a clause requiring that at least edge of the pair is good.  Thus finding an assignment satisfying the axioms is equivalent to solving a 2-satisfiability problem (2 because each clause only involves two edges). Equivalently,
we may consider the graph $H=(E(G), N \cup B)$, where a pair of edges is $\in N$ if they do not touch and a pair of edges is $\in B$ if they induce a path of length 2. We seek a partition of $V(H)=\mathfrak{G}\cup \mathfrak{B}$ such that $\mathfrak{G}$ is a stable set in $(E(H),N)$ and $\mathfrak{B}$ is a stable set in $(E(H),B)$. Such a partition exists if and only if there is a certain kind of alternating walk.  This result is due to Alexander Schrijver \cite{S}.  It is convenient to prove a stronger statement, which we now state.

Let $H_N$ and $H_B$ be graphs on the same vertex set, $V$, with edge sets $N$ and $B$ ($N$ and $B$ need not be disjoint, as they are in the graphs defined above).  A \emph{walk} of length $l$ is a sequence of vertices $v_1, v_2, \ldots, v_{l+1}$, (not necessarily distinct) such that $(v_i, v_{i+1}) \in N\cup B$, $i \in [l]$.
A walk is \emph{closed} if $v_1 = v_{l+1}$. A walk is \emph{alternating} if edges of the form $(v_{2i-1}, v_{2i})$ are in $N$, and edges of the form $(v_{2i}, v_{2i+1})$ are in $B$, (or the same with $N$ and $B$ switched).  Closed alternating walks of odd length are possible with this definition, but if the vertex labels are cyclicly permuted, it is no longer alternating. We will call such a walk (and this name will hold under any cyclic permutation of the vertices) an \emph{AACW} (almost alternating closed walk) with \emph{nose} $v_1$, where $v_1$ is the unique vertex so that the walk can be written as $v_1, v_2, \ldots, v_{l}, v_{l+1}$ and be alternating.  

\begin{theorem}
Exactly one of the following holds:
\setcounter{ctr}{0}
\begin{list}{\emph{(\alph{ctr})}} {\usecounter{ctr} \setlength{\itemsep}{1pt} \setlength{\topsep}{3pt}}
\item There is a partition $V =\mathfrak{G}\cup \mathfrak{B}$ such that $\mathfrak{G}$ is a stable set in $H_N$ and $\mathfrak{B}$ is a stable set in $H_B$.
\item There is an even closed alternating walk such that two vertices an odd distance apart in the walk are identical.
\end{list}
\end{theorem}

\begin{figure}[tbp]
\begin{center}
\includegraphics [scale=1.1]{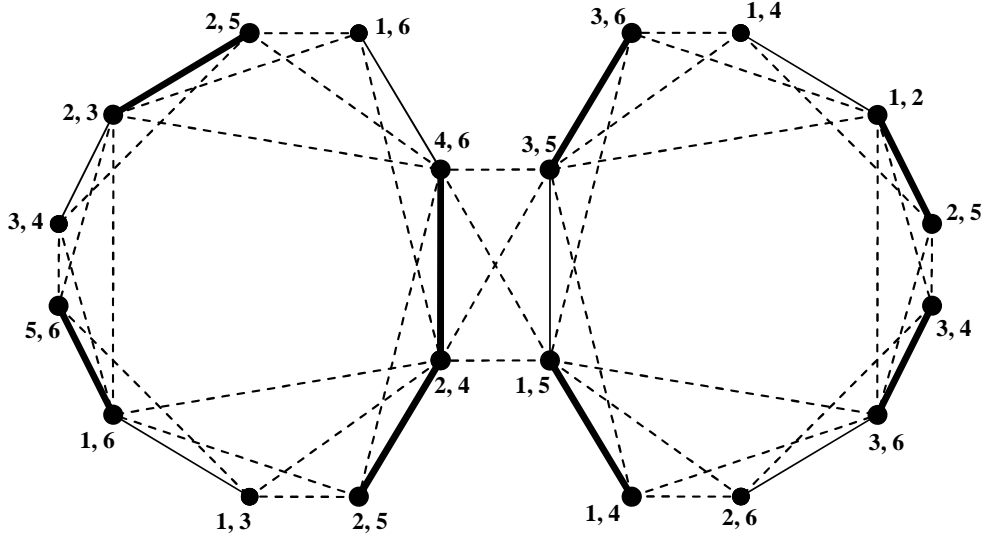}
\end{center}
\caption{A subgraph that makes it impossible to partition $E$ into good and bad edges.  The good and bad edges shown is a failed attempt at an assignment satisfying the good and bad edge axioms.  The two edges on the far right fail to satisfy the good edge axiom.}
\label{2satcountere}
\end{figure}

Figure \ref{2satcountere} is a representation of a walk as described in (b), except this is a drawing of $G$, not $H$.  In $H$, this is two AACW's of length 7; the noses are the two central vertical edges.  An AACW of length 3 or 5 creates an antitriangle in $G$, so this may be the smallest example, but we have not checked carefully.  

At this point we can show that the application of theorem \ref{Goodbad} to graphs with fractional clique covering number less than 3 is in some sense best possible. Suppose there are $k$ cliques of $G \in \mathfrak{A}$, $c_{1},c_{2},\ldots ,c_{k}$ such that every vertex is in at least $\frac{k}{3}$ of them ($G$ has fractional clique covering number at most 3).  If $k$ is odd, we can successfully label edges good and bad.  Label edges using the same rules as before (section 3.1); since there cannot be exactly $\frac{k}{2}$ cliques containing $u$ or $v$, the
same arguments as before show the good and bad edge axioms are satisfied. Observe that if $\overline{G}$ has chromatic number 3, there are three cliques covering $V(G)$, and therefore SSH holds
in this case. 

If $k$ is even, however, $E(G)$ cannot necessarily be
partitioned into good edges and bad edges.  The labels in figure \ref{2satcountere} represent 6 different cliques; every vertex is in exactly 2 of them.  If the only non-edges are those drawn, then this graph is in $\mathfrak{A}$. Evidently, its edges cannot be partitioned into good edges and bad edges satisfying the axioms so theorem \ref{Goodbad} cannot be applied.  Also, it seems that even after deleting dominating edges, we obtain a graph with no large clique or small cutset, although we have not checked this carefully.

\section{Conclusions, conjectures, and future work}

Some questions one might have at this point are ``will the good and bad edge axioms help us say anything about \emph{all} graphs with no antitriangle?,'' ``what happens when the minimum degree is $n-O(n^{4/5})$?,'' ``what about $n$ odd?,'' and ``why the didn't we try induction?''.  We will attempt some answers.

We can only construct random-like graphs with no antitriangle when the average degree is $n-O(n^{1/2})$ or larger.  In this degree range, theorem \ref{n4/5claim} tells us a lot.  Graphs with smaller minimum degree than $n-O(n^{1/2})$ have cliques too large for a typical random graph because the non-neighbors of every vertex are a clique.  This suggests that in this density range, graphs with no antitriangle tend to have structure like that of a smaller c-blown up graph.  It is here that we seek to apply results like theorem \ref{Goodbad} and corollary \ref{goodmatchingcor2}.  So far we are only successful for graphs with fractional clique covering number less than 3, in which case degrees are around $n-\frac{1}{3}n$. So, for example, how might we extend the results using the good and bad edge axioms to graphs with minimum degree $n-O(n^{4/5})$?  Corollary \ref{goodmatchingcor2} gives us one prescription, but what do the resulting medium and bad edges look like?  Perhaps we could compute bounds on the number of pairs of medium edges that do not touch, which might lead to bounds on the size of a complete minor. It is strange that the graphs that give us the most trouble are
denser than the graphs that our results apply to; large complete minors should be easier to find in denser graphs.  In a sense, the problem is not that
it is difficult to find a complete minor in these graphs (ones with minimum degree $n-O(n^{4/5})$, say), but rather that we cannot say anything about them.

Another idea for extending the good and bad edge axioms is to assign weights to the edges.  We may think of the weights as distances.  Edges between vertices with ``similar'' neighbor sets  (those that are close to being c-twins) will receive small weights and will be like bad edges of varying degrees.  Edges that are adjacent to many other edges will receive large weights and will be like good edges of varying degrees.  These two ways of choosing edge weights are similar, but do not agree exactly, and it is not clear what the right weighting function is.  

Let $G$ be as in theorem \ref{goodmatching}. It would be nice if we could modify the proof of theorem \ref{goodmatching} to work for $n$ odd. We think that if $G$ has no dominating edge, small cutset, or large clique, (as in (a), (b), and (c) of theorem \ref{goodmatching}) then we can choose any vertex, $q$, to be a prevertex and use edges for the other prevertices.

Let $Z_q$ be the clique of non-neighbors of $q$. An obvious idea is to
apply Tutte's theorem to the graph $G' \equiv (V-q, \mathfrak{G}-E(Z_q))$. A perfect matching in $G'$ together with $q$ are the prevertices of a complete minor of $G$. A similar proof to that of theorem \ref{goodmatching} works in
quite a few cases, but not all.  We have found a graph with no dominating edge, no small
cutset, and no large clique, and no matching of edges in $\mathfrak{G}-E(Z_q)$ saturates $V-q$. It is possible that if a special set of good and bad edges is chosen, perhaps as in corollary \ref{goodmatchingcor}, then $G'$ does have a perfect matching.  However, our investigations suggest that several natural choices for special sets of good and bad edges do not work.  It seems best to try another approach for $n$ odd.

Another approach to SSH is to suppose there is a cut of size $(n+1)/2$ but no smaller cutset.  Select any vertex $q$.  We have conjectured that there is a minor using only prevertices of size two and the vertex $q$.  Can we show that in this case? Along similar lines, can we show SSH if there is a clique of size $\lceil n/2 \rceil - 1$?  These questions are surprisingly difficult, and solutions may yield insight into the general case.

Tutte's theorem is a fantastic structure theorem for graphs with a perfect matching, but it's not exactly what we need for this problem.  Perhaps we can find an appropriate strengthening of SSH that we can apply inductively to get a perfect matching of medium edges.  We want it to give us a special perfect matching of medium edges, not just any, as Tutte's theorem gives us.  \par
It will require much cleverness to get induction to work on this problem. Suppose we have an edge, $e=(u,v)$, that seems like a good candidate to be a prevertex, and then we inductively obtain the prevertices of a minor on $G \backslash \{u,v\}$.  It is not at all clear that the prevertices of the minor in $G \backslash \{u,v\}$ will touch $e$, so we have to find very special prevertices, not just any.  Labeling some edges bad provides a way for us to exclude \emph{some} sets of prevertices, however, we need something more powerful to tackle the general case.

\section*{Acknowledgements}

I thank Professor Seymour for his many awesome ideas, for teaching me how to do good graph theory, and for scolding me for drawing a planar graph with edge crossings.  I thank Professor Gy\~{o}ri for his help fall semester.  We worked very well together, and I enjoyed struggling with him on this difficult problem.  I thank Professor Sudakov for his helpful discussions about random graphs.

\end{document}